\documentclass{article}

\usepackage{delarray,verbatim,enumerate,a4wide}
\usepackage{amsmath,amsthm,amstext,amsbsy,amssymb,amsfonts,amscd}
\usepackage[all]{xy}
\usepackage[french]{babel}
\usepackage[T1]{fontenc}
\usepackage[utf8]{inputenc}

\usepackage{scalerel}
\usepackage{upgreek,bbm}
\usepackage[small]{titlesec}

\usepackage{color}
\definecolor{vert}{rgb}{0.1,0.4,0.2}
\usepackage[colorlinks=true,linkcolor=blue,citecolor=vert]{hyperref}

\usepackage{calligra}
\DeclareFontShape{T1}{calligra}{m}{n}{<->s*[0.95]callig15}{}
\DeclareMathAlphabet{\mathscr}{T1}{calligra}{m}{n}

\newtheorem{Th}{Théorème}[]
\newtheorem{Lem}[Th]{Lemme}
\newtheorem{Prop}[Th]{Proposition}
\newtheorem{Cor}[Th]{Corollaire}
\newtheorem{Conj}[Th]{Conjecture}
\newtheorem{Sco}[Th]{Scolie}

\newtheorem*{Th*}{Théorème}
\newtheorem*{Lem*}{Lemme}
\newtheorem*{Cor*}{Corollaire}
\newtheorem*{Def*}{Définition}

\def\Preuve{\noindent {\it Preuve.~}}
\def\PreuveTh{\noindent {\it Preuve du Théorème.~}}
\def\Remarque{\smallskip\noindent {\it Remarque.~}}

\def\Nota{\smallskip\noindent {\it Nota.~}}

\def\NN{\mathbb N}	\def\ZZ{\mathbb Z}		
\def\F2{\mathbb{F}_2}	\def\Z2{\mathbb{Z}_2}		
\def\Zl{{\mathbb{Z}_\ell}} 	\def\Ql{{\mathbb{Q}_\ell}}		

\def\I{\mathcal  I} 		\def\P{\mathcal  P}		\def\U{\mathcal  U}	\def\F{\mathcal  F}	
\def\J{\mathcal  J}  		\def\C{\mathcal  C}		\def\R{\mathcal  R}		
 	\def\Pl{\mathcal  P\ell}  	\def\Cl{\mathcal  C\!\ell}	
\def\E{\mathcal  E}		\def\T{\mathcal  T}			\def\D{\mathcal  D}

						\def\G{\mathcal  G}

		\def\p{{\mathfrak p}}				\def\a{{\mathfrak a}}		
		\def\l{{\mathfrak l}}

\def\wi{\widetilde}				
\def\rg{\operatorname{rg}}	
	\def\deg{\operatorname{deg}}		
\def\Gal{\operatorname{Gal}}			
\def\Ker{\operatorname{Ker}}				

\newcommand\scale[2]{\vstretch{#1}{\hstretch{#1}{#2}}}

\newcommand\si[1]{\scale{.7}{#1}}	
\newcommand\ph{{\phantom{*}}}

\def\%{{\scale{.8}{\infty}}}				

\makeatletter
\newcommand*\wt[2][0.2ex]{%
        \begingroup
        \mathchoice{\wt@helper{#1}{#2}{\displaystyle}{\textfont}}
                   {\wt@helper{#1}{#2}{\textstyle}{\textfont}}
                   {\wt@helper{#1}{#2}{\scriptstyle}{\scriptfont}}
                   {\wt@helper{#1}{#2}{\scriptscriptstyle}{\scriptscriptfont}}%
        \endgroup
        #2%
}
\newcommand*\wt@helper[4]{%
        \def\currentfont{\the#41}%
        \def\currentskewchar{\char\the\skewchar\currentfont}%
        \setbox\tw@\hbox{\currentfont$#2$\currentskewchar}%
        \dimen@ii\wd\tw@
        \setbox\tw@\hbox{\currentfont$#2${}\currentskewchar}%
        \advance\dimen@ii-\wd\tw@
        \rlap{\raisebox{-#1}{$\m@th#3\kern\dimen@ii\widetilde{\phantom{#2}}$}}%
}
\makeatother

\def\wE{\,\wt[0.1ex]{\!\mathcal E}}		\def\wU{\wt[0.2ex]{\mathcal U}}	\def\wD{\wt[0.2ex]{\mathcal D}}

\begin{document}

\title{\Large\bf Conjecture cyclotomique et semi-simplicité des modules d'Iwasawa}

\author{ Jean-François {\sc Jaulent} }
\date{}
\maketitle
\bigskip\bigskip\bigskip

{\small
\noindent{\bf Résumé.} Nous montrons que la conjecture cyclotomique sur les modules d'Iwasawa $S$-décomposés $T$-ramifiés, introduite dans un article antérieur et vérifiée par les corps abéliens, gouverne le $\Zl$-rang du sous-module des points fixes et du quotient des genres pour tous   ensembles finis disjoints de places $(S,T)$.\par
Finalement, en cas de conjugaison complexe nous montrons que la forme forte de la conjecture est, tout comme la forme faible, équivalente à la conjonction des conjectures de Leopodt et de Gross-Kuz'min.

\

\noindent{\bf Abstract.} We show that the cyclotomic conjecture on the characteristic polynomial of $T$-ramified $S$-split Iwasawa modules, introduced in a previous paper and satisfied by abelian fields, governs the $\Zl$-rank of the submodule of fixed points for all finite disjoint sets $S$ and $T$ of places.\par
Last, in the CM-case we prove that the weak and the strong versions of the cyclotomic conjecture are both equivalent to the conjunction of the classical conjectures of Leopoldt and Gross-Kuz'min.

\

\noindent{\em Mathematics Subject Classification}: Primary 11R23; Secondary 11R37.\

\noindent{\em Keywords}: cyclotomic conjecture, Iwasawa modules, Leopoldt conjecture, Gross-Kuz'min conjecture
}

\tableofcontents

%%%%%%%%%%%%%%%%%%%%%%%%%%%%%%%%%%%%%%%%%%%%%%%%%%%%%%%%
%%%%%%%%%%%%%%%%%%%%%%%%%%%%%%%%%%%%%%%%%%%%%%%%%%%%%%%%
\section*{Introduction}
\addcontentsline{toc}{section}{Introduction}
\medskip
%%%%%%%%%%%%%%%%%%%%%%%%%%%%%%%%%%%%%%%%%%%%%%%%%%%%%%%%
%%%%%%%%%%%%%%%%%%%%%%%%%%%%%%%%%%%%%%%%%%%%%%%%%%%%%%%%

%%%%%%%%%%%%%%%%%%%%%%%%%%%%%%%%%%%%%%%%%%%%%%%%%%%%%%%%%

Pour chaque nombre premier $\ell$, il est avancé dans  \cite{J63} la conjecture générale suivante:

\begin{Conj}[\bf Conjecture cyclotomique forte]
Soient $K$ un corps de nombres arbitraire; $K_{\si{\infty}}=\bigcup_{n\in\NN}K_n$ sa $\Zl$-extension cyclotomique; $\Lambda$ l'algèbre d'Iwasawa  de $\Gamma=\Gal(K_{\si{\infty}}/K)=\gamma^{\Zl}$; et $Pl^\ell_{\si{K}}=S_{\si{K}} \sqcup T_{\si{K}}$ une partition de  l'ensemble $Pl^\ell_{\si{K}}$ des places de $K$ au-dessus de $\ell$.\par

Alors le polynôme caractéristique du $\Lambda$-module $\,\C^{T_{\si{K}}}_{S_{\si{K}}}(K_{\si{\infty}})=\Gal(H^{T_{\si{K}}}_{S_{\si{K}}}(K_{\si{\infty}})/K_{\si{\infty}})$ attaché à la pro-$\ell$-extension abélienne  maximale $H^{T_{\si{K}}}_{S_{\si{K}}}(K_{\si{\infty}})$ de $K_{\si{\infty}}$ qui est  $S_{\si{K}}$-décomposée et $T_{\si{K}}$-ramifiée n'est pas divisible par $\omega=\gamma-1$.
 Autrement dit, son sous-module des points fixes  $\,\C^{T_{\si{K}}}_{S_{\si{K}}}(K_{\si{\infty}})^\Gamma$ est fini:\smallskip
 
 \centerline{$\C^{T_{\si{K}}}_{S_{\si{K}}}(K_{\si{\infty}})^\Gamma \sim 1$.}
\end{Conj}

Cette conjecture cyclotomique contient celles de Leopoldt et de Gross-Kuz'min (cf. e.g. \cite{J55}), qui correspondent respectivement aux cas $(S_{\si{K}},T_{\si{K}})=(\emptyset,Pl_{\si{K}}^\ell)$ et $(S_{\si{K}},T_{\si{K}})=(Pl_{\si{K}}^\ell,\emptyset)$ et réciproquement: si $K$ est totalement réel, ou encore si $K$ est un corps à conjugaison complexe extension quadratique totalement imaginaire d'un sous-corps totalement réel, la conjecture cyclotomique ci-dessus est vraie dès lors que les ensembles $S_{\si{K}}$ et $T_{\si{K}}$ sont stables par conjugaison complexe et que $K$ vérifie à la fois la conjecture de Leopoldt et celle de Gross-Kuz'min (\cite{J63}, Th. 5). En résumé, 
%%%%%%%%%%%%%%%%%%%%%%%%%%%%%%%%%%%%%%%%%%%%%%%%%
{\em restreinte aux seules partitions $Pl^\ell_{\si{K}}=S_{\si{K}} \sqcup T_{\si{K}}$ stables par conjugaison complexe},
%%%%%%%%%%%%%%%%%%%%%%%%%%%%%%%%%%%%%%%%%%%%%%%%%
elle est %alors 
équivalente à la conjonction des conjectures de Leopoldt et de Gross-Kuz'min.

Pour prendre en compte cette restriction sur la stabilité par conjugaison des ensembles $S_{\si{K}}$ et $T_{\si{K}}$, restée ambiguë dans \cite{J63}, introduisons la forme affaiblie ci-après de la conjecture:

\begin{Conj}[\bf Conjecture cyclotomique faible]
Nous disons qu'un corps de nombres $K$ à conjugaison complexe (i.e. extension quadratique totalement imaginaire d'un sous-corps totalement réel $K^+$) satisfait la conjecture cyclotomique faible pour un nombre premier $\ell$, lorsqu'il vérifie la conjecture cyclotomique pour toute partition  $Pl^\ell_{\si{K}}=S_{\si{K}} \sqcup T_{\si{K}}$ stable par la conjugaison complexe.
\end{Conj}

Avec ces conventions le Théorème 5 de \cite{J63} affirme alors qu'un tel $K$ vérifie la Conjecture cyclotomique {\em faible} pour $\ell$ si et seulement s'il vérifie simultanément la conjecture de Leopoldt et celle de Gross-Kuz'min pour ce même $\ell$; ce qui est toujours le cas pour $K$ abélien.
\medskip

Le but de la présente note est de déterminer le $\Zl$-rang du module des points fixes $\,\C^{T_{\si{K}}}_{S_{\si{K}}}(K_\%)^\Gamma$, 
sous l'une ou l'autre des conjectures cyclotomiques ci-dessus, pour tout couple $(S_{\si{K}},T_{\si{K}})$ d'ensembles finis disjoints de places de $K$, lors même que la condition $Pl_{\si{K}}^\ell\subset S_{\si{K}}\cup T_{\si{K}}$ n'est pas satisfaite, et d'en tirer quelques conséquences sur les modules d'Iwasawa $\,\C^{T_{\si{K}}}_{S_{\si{K}}}(K_\%)$ pour $K$ à conjugaison complexe en liaison avec les résultats de \cite{J43,JMa,JMP}. Nous procédons pour cela en trois temps:\smallskip

-- Dans le cadre de la conjecture faible d'abord, qui permet de traiter le cas où $S_{\si{K}}$ et $T_{\si{K}}$ sont stables par la conjugaison complexe.
La formule de rang que nous obtenons (Th. \ref{TP1}), qui est donc vérifiée dans le cas abélien, peut être regardée comme une généralisation de la classique conjecture de Coates et Lichtenbaum démontrée par Greenberg dans ce même contexte (cf. \cite{CL,Grb1}).\smallskip

-- Dans le cadre de la conjecture forte ensuite qui ouvre sur le cas plus général où n'est pas fait d'hypothèse de stabilité. Pour cela, nous commençons par tirer quelques conséquences algébriques de la Conjecture cyclotomique généralisant le résultat de semi-simplicité de Greenberg (cf. \cite{Grb1}) qui prouve que, pour un corps abélien $K$, le polynôme minimal du $\Lambda$-module $\Gal(H_{Pl_{\si{K}}^\ell}(K_\%)/K_{\si{\infty}})$ de la plus grande pro-$\ell$-extension abélienne non ramifiée et $\ell$-décomposée $H_{Pl_{\si{K}}^\ell}(K_\%)$ de $K_{\si{\infty}}$,  n'est pas divisible par $(\gamma-1)$. Il en résulte une formule de rang qui étend la précédente (Th. \ref{TP2}). \smallskip

-- Cela fait, nous abordons la comparaison des formes faible et forte de la conjecture et montrons qu'elles sont en fait équivalentes (Th. \ref{TP3}).
L'idée centrale est la suivante: on s'intéresse à un $\Zl$-module qui n'est pas {\em a priori} stable par conjugaison complexe, mais qui est canoniquement l'image d'un autre module stable par conjugaison, dont on peut donc définir les composantes réelle et imaginaire. Sous la conjecture de Leopoldt, il se trouve que sa composante réelle est pseudo-nulle. Quotientant donc le module de départ par l'image de celle-ci, on obtient un module pseudo-isomorphe au module de départ lequel est, lui, purement imaginaire.\smallskip

Enfin ans la dernière partie de cette note nous abordons directement le calcul du $\Zl$-rang à l'aide de la théorie $\ell$-adique du corps de classes. La preuve alternative particulièrement concise que nous en donnons (Cor. \ref{CF}) redonne immédiatement l'équivalence établie plus haut.
\smallskip

Signalons pour finir qu'au cours de l'élaboration de ce résultat nous avons eu connaissance d'un travail indépendant de Lee et Yu \cite{LY} complétant  l'étude antérieure de Lee et Seo \cite{LS} sur une forme équivalente de la conjecture cyclotomique. Avec celles données dans la présente note, on dispose de ce fait de trois façons différentes de déterminer le $\Zl$-rang du module $\C^{T_{\si{K}}}_{S_{\si{K}}}(K_{\si{\infty}})^\Gamma$, en présence d'une conjugaison complexe sous les conjectures équivalentes ci-dessus.\medskip

\Remarque Comme expliqué dans \cite{J63}, Scolie 8, les places étrangères à $\ell$ étant sans incidence sur le $\Zl$-rang de $\,\C^{T_{\si{K}}}_{S_{\si{K}}}(K_\%)^\Gamma$, il est toujours possible de supposer $S_{\si{K}} \cup T_{\si{K}} \subset Pl_{\si{K}}^\ell$ dans les démonstrations, sans restreindre aucunement la généralité.

%%%%%%%%%%%%%%%%%%%%%%%%%%%%%%%%%%%%%%%%%%%%%%%%%%%%%%%%%%%%
%%%%%%%%%%%%%%%%%%%%%%%%%%%%%%%%%%%%%%%%%%%%%%%%%%%%%%%%

\newpage
%%%%%%%%%%%%%%%%%%%%%%%%%%%%%%%%%%%%%%%%%%%%%%%%%%%%%%%%%%%
\section{Théorème Principal sous la conjecture faible}
%%%%%%%%%%%%%%%%%%%%%%%%%%%%%%%%%%%%%%%%%%%%%%%%%%%%%%%%%%%

Supposons maintenant que $K$ soit une extension quadratique totalement imaginaire d'un sous-corps totalement réel $K^+$. Notons $\tau$ la conjugaison complexe et $\Delta=\Gal(K/K^+)=\{1,\tau\}$.\smallskip

Rappelons que, pour $\ell$ impair, chaque $\Zl[\Delta]$-module $M$ s'écrit comme somme directe de ses composantes réelle $M^+=M^{e_{\si{+}}}$ et imaginaire $M^-=M^{e_{\si{-}}}$ via les idempotents  $e_\pm= \frac{1}{2}(1 \pm \tau)$.  Et, si $\ell$ vaut $2$ et si $M$ est $\Z2$-noethérien, confondre le noyau de $(1\pm\tau)$ avec l'image de $(1\mp\tau)$ donne lieu à une erreur finie.
On peut donc dans tous les cas définir composantes réelle et imaginaire de $M$ comme image et noyau respectifs de $(1+\tau)$ et écrire à un fini près: $M \sim M^+\oplus M^-$.

Pour chaque ensemble de places $S$ de $K^+$, notons de même $S^-$ le sous-ensemble de celles qui sont décomposées par la conjugaison complexe et $S^+$ son complémentaire.\smallskip

Avec ces conventions, sous la conjecture faible le {\em Théorème principal} s'énonce comme suit:

\begin{Th}\label{TP1}
Soient $\ell$ un nombre premier; $K$ une extension quadratique totalement imaginaire d'un sous-corps $K^+$ totalement réel; $S_{\si{K^+}}$ et $T_{\si{K^+}}$ deux ensembles finis disjoints de places de $K^+$; $R_{\si{K^+}}=Pl_{\si{K^+}}^\ell\setminus (S_{\si{K^+}}\cup T_{\si{K^+}})$ l'ensemble des places au-dessus de $\ell$ qui ne sont ni dans $S_{\si{K^+}}$ ni dans $T_{\si{K^+}}$.\par
Soient enfin $\,\C^T_S(K_{\si{\infty}})=\Gal(H^T_S(K_{\si{\infty}})/K_{\si{\infty}})$ le groupe de Galois attaché à la pro-$\ell$-extension abélienne  maximale $H^T_S(K_{\si{\infty}})$ de $K_{\si{\infty}}$ qui est  $S$-décomposée et $T$-ramifiée; et $\,\C^T_S(K_{\si{\infty}})^\Gamma$ le sous-module de $\,\C^T_S(K_{\si{\infty}})$ fixé par $\Gamma=\Gal(K_\%/K)$.
Sous la conjecture cyclotomique faible, il vient:\smallskip
\begin{itemize}
\item[(i)] La composante réelle du groupe $\,\C^T_S(K_{\si{\infty}})^\Gamma$ est finie: $\big(\C^T_S(K_{\si{\infty}})^\Gamma\big)^+\sim\, \C^T_S(K^+_{\si{\infty}})^\Gamma \sim 1$.
\item[(ii)] Sa composante imaginaire est un $\Zl$-module de rang: $\rg_{\Zl}\big(\C^T_S(K_{\si{\infty}})^\Gamma\big)^-=|R^-|$,
où $|R^-|$ désigne le nombre de places de $R$ qui sont décomposées par la conjugaison complexe.
\end{itemize}
\end{Th}

\Preuve Comme rappelé plus haut, nous pouvons supposer, sans perte de généralité, que $R \sqcup S \sqcup T$ et une partition de $Pl^\ell$. De plus, côté réel, nous avons immédiatement: $\big(\C^T_S(K_{\si{\infty}})^\Gamma\big)^+\sim\, \C^T_S(K^+_{\si{\infty}})^\Gamma$, puisque l'opérateur $1+\tau$ correspond à la norme $N_{\si{K/K^+}}$.\smallskip

Considérons alors le quotient des genres $^\Gamma\C^T_S(K^+_{\si{\infty}})=\Gal(H^T_S(K^+_\%/K^+)/K_\%^+)$, où $H^T_S(K^+_\%/K^+)$ désigne la plus grande sous-extension de $H^T_S(K^+_\%)$ qui est abélienne sur $K^+$. Par construction,  $H^T_S(K^+_\%/K^+)$ est $\ell$-ramifiée sur $K_\%^+$, donc sur $K^+$ et, par conséquent, de degré fini sur $K_\%^+$ sous la conjecture de Leopoldt (qui est vérifiée ici, puisque la conjecture cyclotomique est supposée l'être).
En résumé, il vient: $^\Gamma\C^T_S(K^+_{\si{\infty}}) \sim 1$; donc, a fortiori: $\,\C^T_S(K^+_{\si{\infty}})^\Gamma \sim 1$.
\smallskip

Regardons maintenant la composante imaginaire. D'après la Proposition 2 de \cite{J63} appliquée aux étages finis $K_n/K$ de la $\Zl$-extension cyclotomique de $K$, l'isomorphisme de modules galoisiens\smallskip

\centerline{$\Cl^T_S(K_n)^\Gamma/cl^T_S\big(\D^T_S(K_n)^\Gamma\big) \simeq \big(\E^T_S(K)\cap N_{\si{K_n/K}}(\R(K_n))\big)/N_{\si{K_n/K}}(\E_S^T(K_n))$}\smallskip

\noindent identifie le quotient du pro-$\ell$-groupe des $S$-classes $T$-infinitésimales ambiges de $K_n$ par le sous-groupe des classes des $S$-diviseurs $T$-infinitésimaux ambiges à un certain quotient du groupe des $S$-unités $T$-infinitésimales; ce qui donne, par passage à la limite projective pour la norme:\smallskip

\centerline{$\C^T_S(K_\%)^\Gamma/\varprojlim cl^T_S\big(\D^T_S(K_n)^\Gamma\big) \simeq \big(\E^T_S(K)\cap \wE(K)\big)/\big(\bigcap_{n\in\NN}N_{\si{K_n/K}}(\E_S^T(K_n)\big)$.}\smallskip

Le point essentiel ici est que, sous la conjecture de Gross-Kuz'min (donc ici encore sous la conjecture cyclotomique), la composante imaginaire du groupe $\,\wE(K)=\bigcap_{n\in\NN} N_{\si{K_n/K}}(\R(K_n))$ des unités logarithmiques se réduit au $\ell$-sous-groupe des racines de l'unité (cf. e.g. \cite{J28}). Il suit:\smallskip

\centerline{$\big(\C^T_S(K_\%)^\Gamma\big)^-\!\sim \varprojlim\, cl^T_S \big(\D^T_S(K_n)^\Gamma\big)^-$.}\smallskip

\noindent Or, le pro-$\ell$-groupe $\big(\D^T_S(K_n)^\Gamma\big)^-$ des $S$-diviseurs étrangers à $T$, qui sont ambiges et imaginaires, est engendré par le sous-groupe $\D^T_S(K)^-$  étendu de $K$ (donc sans incidence sur la limite projective) et le sous-groupe $\D^T_S(K_n)^{[R]-}$ construit sur les produits $\a_n(\p)=\prod_{\p_n\mid \p}\p_n$, pour $\p\in R^-$, lesquels satisfont les identités normiques $N_{\si{K_m/K_n}} (\a_m(\p))=\a_n(\p)$, pour $m\ge n\gg 0$.
Son sous-groupe principal  $\P^T_S(K_n)^{[R]-}$ étant ultimement constant, puisque les $S$-unités imaginaires $T$-infinitésimales des $K_n$ proviennent d'un $K_{n_{\si{0}}}$ aux racines de l'unité près, il vient bien finalement:\smallskip

\centerline{$\big(\C^T_S(K_\%)^\Gamma\big)^-\!\sim \varprojlim\, cl^T_S(\D^T_S(K_n)^{[R]-}) \sim \varprojlim\, \D^T_S(K_n)^{[R]-} \simeq  \Zl^{|R^-|}$.}

\newpage
%%%%%%%%%%%%%%%%%%%%%%%%%%%%%%%%%%%%%%%%%%%%%%%%%%%%%%%%%%%
\section{Semi-simplicité des modules d'Iwasawa $\,\C^T_S(K_\infty)$}
%%%%%%%%%%%%%%%%%%%%%%%%%%%%%%%%%%%%%%%%%%%%%%%%%%%%%%%%%%%

Avant d'énoncer le Théorème Principal dans le cadre plus général de la conjecture forte, commençons par préciser quelques conséquences algébriques de celle-ci.\smallskip

Notons $\gamma$ un générateur topologique du groupe de Galois $\Gamma=\Gal(K_\%/K)$ et $\Lambda=\Zl[[\gamma-1]]$ l'algèbre d'Iwasawa attachée à $\Gamma$. Le point essentiel est que les facteurs cyclotomiques du polynôme minimal $\Pi^{T_{\si{K}}}_{S_{\si{K}}}(\gamma-1)$ du sous-module de $\Lambda$-torsion $\,\T^{T_{\si{K}}}_{S_{\si{K}}}(K_\%)$ de $\C^{T_{\si{K}}}_{S_{\si{K}}}(K_\%)$ sont de multiplicité 1:

\begin{Prop}\label{Prop}
Etant donné un corps de nombres $K$  qui satisfait la conjecture cyclotomique forte pour $\ell$, soient $K_\%=\bigcup K_n$ sa $\Zl$-extension cyclotomique; $S_{\si{K}}$ et $T_{\si{K}}$ deux ensembles finis disjoints quelconques de places finies de $K$; puis $R_{\si{K}}= Pl_{\si{K}}^\ell\setminus(S_{\si{K}} \cup T_{\si{K}})$ le sous-ensemble des places de $K$ au-dessus de $\ell$ qui ne sont ni dans $S_{\si{K}}$ ni dans $T_{\si{K}}$.
Soient enfin $\,\C^{T_{\si{K}}}_{S_{\si{K}}}=\C^{T_{\si{K}}}_{S_{\si{K}}}(K_\%)$ le groupe de Galois $\Gal(H^{T_{\si{K}}}_{S_{\si{K}}}(K_\%)/K_\%)$ de la pro-$\ell$-extension abélienne $S_{\si{K}}$-décomposée $T_{\si{K}}$-ramifiée maximale de $K_\%$ et $\,\C^{T_{\si{K}}}_{S_{\si{K}}}[R_{\si{K}}]$ engendré par les sous-groupes de décomposition des places de $R_{\si{K}}$.
\smallskip
\begin{itemize}
\item[(i)] Le sous-module de $\,\C^{T_{\si{K}}}_{S_{\si{K}}}$ fixé par $\Gamma=\Gal(K_\%/K)$ contient $\,\C^{T_{\si{K}}}_{S_{\si{K}}}[R_{\si{K}}]^\Gamma$ avec un indice fini:\smallskip

\centerline{$\C^{T_{\si{K}}}_{S_{\si{K}}}[R_{\si{K}}]^\Gamma \sim \;(\C^{T_{\si{K}}}_{S_{\si{K}}})^\Gamma$.}\smallskip

\item[(ii)] Et celui du quotient $\,\C^{T_{\si{K}}}_{R_{\si{K}}S_{\si{K}}}=\C^{T_{\si{K}}}_{S_{\si{K}}}/\C^{T_{\si{K}}}_{S_{\si{K}}}[R_{\si{K}}]$ est fini: $(\C^{T_{\si{K}}}_{R_{\si{K}}S_{\si{K}}})^\Gamma\sim 1$.
\end{itemize}\smallskip

\noindent En particulier, le polynôme minimal $\Pi^{T_{\si{K}}}_{S_{\si{K}}}(\gamma-1)$ du sous-module de $\Lambda$-torsion $\,\T^{T_{\si{K}}}_{S_{\si{K}}}$ de $\,\C^{T_{\si{K}}}_{S_{\si{K}}}$ n'est pas divisible par $\Phi_1(\gamma)^2=(\gamma-1)^2$. Plus généralement, sous la conjecture cyclotomique dans $K_\%$, le polynôme $\Pi^{T_{\si{K}}}_{S_{\si{K}}}(\gamma-1)$ n'est divisible par aucun carré de la forme $\Phi_{\ell^n}(\gamma)^2=((\gamma^{\ell^n}-1)/(\gamma^{\ell^{n-\si{1}}}-1))^2$.
\end{Prop}

\Preuve Prenant les points fixes par $\Gamma$ dans la suite exacte courte qui définit $\,\C^{T_{\si{K}}}_{S_{\si{K}}}[R_{\si{K}}]$, \smallskip

\centerline{$1 \to \C^{T_{\si{K}}}_{S_{\si{K}}}[R_{\si{K}}] \to \C^{T_{\si{K}}}_{S_{\si{K}}} \to \C^{T_{\si{K}}}_{R_{\si{K}}S_{\si{K}}} \to 1$,}

\noindent nous obtenons la suite:

\centerline{$1 \to( \C^{T_{\si{K}}}_{S_{\si{K}}}[R_{\si{K}}])^\Gamma \to (\C^{T_{\si{K}}}_{S_{\si{K}}})^\Gamma \to (\C^{T_{\si{K}}}_{R_{\si{K}}S_{\si{K}}})^\Gamma$;}\smallskip

\noindent et la conjecture cyclotomique appliquée avec $R_{\si{K}}\cup S_{\si{K}}$ et $T_{\si{K}}$ nous donne: $(\C^{T_{\si{K}}}_{R_{\si{K}}S_{\si{K}}})^\Gamma \sim 1$.
Les deux assertions $(i)$ et $(ii)$ en résultent immédiatement.\smallskip

Il suit de là que $(\gamma-1)$ et $(\gamma-1)^2$ ont même noyau dans $\,\T^{T_{\si{K}}}_{S_{\si{K}}}$, i.e. que $\Phi_1(\gamma)=\gamma-1$ apparait dans $\Pi^{T_{\si{K}}}_{S_{\si{K}}}(\gamma-1)$ avec une multiplicité au plus 1; et, sous la conjecture cyclotomique dans $K_n$, qu'il en va de même des $\Phi_{\ell^m}(\gamma)$ pour $m\le n$.

\begin{Cor}
Sous la conjecture forte le sous-module des points fixes $(\C^{T_{\si{K}}}_{S_{\si{K}}})^\Gamma=(\T^{T_{\si{K}}}_{S_{\si{K}}})^\Gamma$ de $\,\C^{T_{\si{K}}}_{S_{\si{K}}}$ est un pseudo-facteur direct du sous-module de $\Lambda$-torsion $\,\T^{T_{\si{K}}}_{S_{\si{K}}}$ de $\,\C^{T_{\si{K}}}_{S_{\si{K}}}$.
\end{Cor}

\Preuve Le quotient $\bar\Pi^{T_{\si{K}}}_{S_{\si{K}}}(\gamma-1)=\Pi^{T_{\si{K}}}_{S_{\si{K}}}(\gamma-1)/(\gamma-1)$ étant étranger à $\gamma-1$ dans l'anneau $\Ql[\gamma-1]$, le produit   de leurs noyaux respectifs dans $\,\C^{T_{\si{K}}}_{S_{\si{K}}}$ est d'indice fini dans $\,\T^{T_{\si{K}}}_{S_{\si{K}}}$ et pseudo-direct (en ce sens que leur intersection $\Ker \bar\Pi^T_S(\gamma-1) \cap\Ker (\gamma-1)$ est finie): $\Ker \bar\Pi^{T_{\si{K}}}_{S_{\si{K}}}(\gamma-1) \times \Ker (\gamma-1) \sim \,\T^T_S$.

\begin{Cor}\label{C}
Soient $S'_{\si{K}}\subset S_{\si{K}}$ et $T'_{\si{K}}\supset T_{\si{K}}$ deux autres ensembles finis disjoints de places de $K$.
Sous la conjecture cyclotomique forte, si $\,\C^{T_{\si{K}}}_{S_{\si{K}}}$ et $\,\C^{T'_{\si{K}}}_{S'_{\si{K}}}$ ont même $\Lambda$-rang,  la surjection canonique $f$ de $\,\C^{T'_{\si{K}}}_{S'_{\si{K}}}$ sur $\,\C^{T_{\si{K}}}_{S_{\si{K}}}$ induit un pseudo-épimorphisme de $(\C_{S'_{\si{K}}}^{T'_{\si{K}}})^\Gamma $vers $(\C^{T_{\si{K}}}_{S_{\si{K}}})^\Gamma$.
\end{Cor}

\Preuve Notons $\,\T^T_S$ le sous-module de $\Lambda$-torsion de $\,\C^T_S$ et $\,\T^{T'}_{S'}$ celui de $\,\C^{T'}_{S'}$. La surjection canonique $f$ envoie $(\C^{T'}_{S'})^\Gamma=(\T^{T'}_{S'})^\Gamma$ vers $(\C^{T}_{S})^\Gamma=(\T^{T}_{S})^\Gamma$.
Et l'identité des $\Lambda$-rangs $\rg_\Lambda\C^{T'}_{S'}=\rg_\Lambda\C^{T}_{S}$ nous assure que $f$ envoie $\,\T^{T'}_{S'}$ {\em sur} $\,\T^{T}_{S}$.
En particulier, le polynôme minimal $\Pi(\gamma-1)$ de $\,\T^{T}_{S}$ divise donc le polynôme minimal $\Pi'(\gamma-1)$ de $\,\T^{T'}_{S'}$. Cela étant:\smallskip
\begin{itemize}
\item Si $\gamma-1$ ne divise pas $\Pi'(\gamma-1)$, le sous-module $(\C^{T}_{S})^\Gamma$ est fini; et il n'y a rien à démontrer.\smallskip
\item Sinon, soit $\bar\Pi'(\gamma-1)=\Pi'(\gamma-1)/(\gamma-1)$, avec $\bar\Pi'(\gamma-1)$ et $\gamma-1$ copremiers. Il vient:\smallskip

\centerline{$(\C^{T}_{S})^\Gamma=(\T^{T}_{S})^\Gamma \sim (\T^{T}_{S})^{\bar\Pi'(\gamma-1)}=f((\T^{T'}_{S'})^{\bar\Pi'(\gamma}) \sim f((\T^{T'}_{S'})^\Gamma) = f((\C^{T'}_{S'})^\Gamma)$.}

\end{itemize}

\newpage
%%%%%%%%%%%%%%%%%%%%%%%%%%%%%%%%%%%%%%%%%%%%%%%%%%%%%%%%%%%%%%%%%%%%%%%
\section{Théorème principal sous la conjecture forte}
%%%%%%%%%%%%%%%%%%%%%%%%%%%%%%%%%%%%%%%%%%%%%%%%%%%%%%%%%%%%%%%%%%%%%%%

Dans le Théorème \ref{TP1}, les ensembles respectifs de places de $K$ au-dessus de $S$ et $T$ sont, du fait même de leur construction, stables par la conjugaison complexe $\tau$. Mais il est facile de s'affranchir de cette restriction, ce qui donne le résultat de Lee et Yu obtenu (pour $\ell$ impair) dans \cite{LY}:

\begin{Th}\label{TP2}
Soient $\ell$ un nombre premier et $K$ une extension quadratique totalement imaginaire d'un sous-corps $K^+$ totalement réel, supposée satisfaire la Conjecture cyclotomique forte pour $\ell$. Étant donnés deux ensembles finis disjoints $S_{\si{K}}$ et $T_{\si{K}}$ de places de $K$, notons $\hat S_{\si{K}}=S_{\si{K}} \cup S_{\si{K}}^\tau$ et $\check T_{\si{K}}=T_{\si{K}} \cap\, T^\tau_{\si{K}}$ leurs saturés respectivement supérieur et inférieur pour la conjugaison complexe $\tau$; désignons par $\hat S_{\si{K^+}}$ et $\check T_{\si{K^+}}$ les ensembles de places de $K^+$ au-dessous de $\hat S_{\si{K}}$ et  $\check T_{\si{K}}$; et notons enfin $\bar R_{\si{K^+}}=Pl_{\si{K^+}}^\ell\setminus (\hat S_{\si{K^+}}\cup \check T_{\si{K^+}})$ l'ensemble des places de $K^+$ au-dessus de $\ell$ qui ne sont ni dans $\hat S_{\si{K^+}}$ ni dans $\check T_{\si{K^+}}$.
Alors le $\Zl$-rang du sous-module ambige $\,\C^{T_{\si{K}}}_{S_{\si{K}}}(K_{\si{\infty}})^\Gamma$ du groupe de Galois $\,\C^{T_{\si{K}}}_{S_{\si{K}}}(K_{\si{\infty}})$ attaché à la pro-$\ell$-extension abélienne $H^{T_{\si{K}}}_{S_{\si{K}}}(K_{\si{\infty}})$  $S_{\si{K}}$-décomposée $T_{\si{K}}$-ramifiée maximale sur $K_{\si{\infty}}$ est:\smallskip

\centerline{$\rg_{\Zl}\,\C^{T_{\si{K}}}_{S_{\si{K}}}(K_{\si{\infty}})^\Gamma =
\rg_{\Zl}\,\C^{\check T_{\si{K}}}_{\hat S_{\si{K}}}(K_{\si{\infty}})^\Gamma=|\bar R_{\si{K^+}}^-|$,}\smallskip

\noindent où $|\bar R_{\si{K^+}}^-|$ désigne le nombre de places de $\bar R_{\si{K^+}}$ qui sont décomposées par la conjugaison complexe. 
\end{Th}

Pour établir ce dernier résultat, nous allons nous appuyer sur l'égalité des rangs:

\begin{Lem}\label{L}
Les $\Lambda$-modules $\,\C^{T_{\si{K}}}_{S_{\si{K}}}(K_{\si{\infty}})$ et $\,\C^{\check T_{\si{K}}}_{\hat S_{\si{K}}}(K_{\si{\infty}})$ ont inconditionnellement même $\Lambda$-rang:\smallskip

\centerline{$\rho^{\si{T}}_{\si{S}} = \rg_\Lambda\,\C^{T_{\si{K}}}_{S_{\si{K}}}(K_{\si{\infty}})  = \rg_\Lambda\,\C^{\check T_{\si{K}}}_{\hat S_{\si{K}}}(K_{\si{\infty}}) = \deg_\ell \check T_{\si{K}}$,}\smallskip

\noindent où $\deg_\ell \check T_{\si{K}} = \sum_{\l\in(\Pl^\ell_{\si{K}} \cap \check T_{\si{K}})} [K_\l^+:\Ql]$ désigne le degré en $\ell$ de l'ensemble de places $\check T_{\si{K}}$.
\end{Lem}

\Preuve Si $K$ contient une racine $\ell$-ième primitive de l'unité $\zeta_\ell$, c'est le Théorème 9 de \cite{JMa}. Sinon, écrivant $K=K^+[\sqrt\delta]$ avec $\delta$ totalement négatif, on peut remplacer $K^+$ par $K^+[\zeta_\ell+\bar\zeta_\ell, \sqrt\delta(\zeta_\ell-\bar\zeta_\ell)]$ et $K$ par $K'=K[\zeta_\ell]$; appliquer le Théorème 2.7 de \cite{J43} à $K'$; et redescendre le résultat dans $K$.

\medskip

\PreuveTh Comme précédemment, nous pouvons supposer sans perte de généralité $S_{\si{K}}$ et $T_{\si{K}}$ contenus dans l'ensemble $Pl^\ell_{\si{K}}$ des places de $K$ au-dessus de $\ell$. Cela étant, les sous-ensembles\smallskip

\centerline{$S_{\si{K}}=\check S_{\si{K}} \sqcup (T_{\si{K}}\cap S_{\si{K}}^\tau) \sqcup S^\circ_{\si{K}};
\quad  T_{\si{K}}=\check T_{\si{K}} \sqcup (S_{\si{K}}\cap T_{\si{K}}^\tau) \sqcup T^\circ_{\si{K}};
\quad R_{\si{K}}=\check R_{\si{K}} \sqcup S_{\si{K}}^{\circ\,\tau} \sqcup T_{\si{K}}^{\circ\,\tau}$.}\smallskip

\noindent forment alors une partition de $Pl_{\si{K}}^\ell$. Il vient: $\bar R_{\si{K}}= \check R_{\si{K}} \sqcup (T^\circ_{\si{K}} \sqcup T^{\circ\,\tau}_{\si{K}})$ puis $\bar R_{\si{K}}^-= \check R_{\si{K}}^- \sqcup (T^\circ_{\si{K}} \sqcup T^{\circ\,\tau}_{\si{K}})$.\smallskip

Procédons par minoration et majoration.\smallskip

$(i)$ Le Corollaire \ref{C}, nous donne directement la minoration: $\rg_\Zl\,\C_{S_{\si{K}}}^{T_{\si{K}}}(K_\%)^\Gamma \ge \rg_\Zl\,\C_{\hat S_{\si{K}}}^{\check T_{\si{K}}}(K_\%)^\Gamma$. Et ce dernier est donné par le Théorème \ref{TP1}: $\rg_\Zl\,\C_{\hat S_{\si{K}}}^{\check T_{\si{K}}}(K_\%)^\Gamma  = |(\bar R_{\si{K^+}}^-|$.\smallskip

$(ii)$ D'autre part, par la Proposition \ref{Prop}, les classes invariantes de $\,\C_{S_{\si{K}}}^{T_{\si{K}}}(K_\%)^\Gamma$ proviennent des familles projectives $\a_n(\p)$ pour $\p\in R_{\si{K}}=\check R_{\si{K}} \sqcup S_{\si{K}}^{\circ\,\tau} \sqcup T_{\si{K}}^{\circ\,\tau}$. Or, par le Théorème \ref{TP1}:
\begin{itemize}
\item pour $\p \in \check R_{\si{K}}$, on a: $\a(\p)\a(\p^\tau)=\a(\p\p^\tau)\sim 1$ dans $\,\C^{\hat T_{\si{K}}}_\ph(K_\%)$ donc, a fortiori, dans $\,\C^{T_{\si{K}}}_{S_{\si{K}}}(K_\%)$;
\item pour $\p \in S_{\si{K}}^{\circ\,\tau}$, on a de même: $\a(\p)\a(\p^\tau)=\a(\p\p^\tau)\sim 1$ dans $\,\C^{\hat T_{\si{K}}}_\ph(K_\%)$ donc dans $\,\C^{T_{\si{K}}}_\ph (K_\%)$; et finalement: $\a(\p)\sim\a(\p^{-\tau})\sim 1$ dans $\,\C^{T_{\si{K}}}_{S_{\si{K}}}(K_\%)$.
\end{itemize}\smallskip

Il vient donc: $\rg_\Zl\,\C_{S_{\si{K}}}^{T_{\si{K}}}(K_\%)^\Gamma\le \frac{1}{2}|\check R_{\si{K}}^-|+|T_{\si{K}}^{\circ\,\tau}|= |\bar R_{\si{K^+}}^-|$.\smallskip

En fin de compte, les deux inégalités réunies nous donnent l'égalité attendue.

\begin{Cor}
Sous les hypothèses du Théorème, le $\Zl$-rang du quotient des genres %$^\Gamma \C^{T_{\si{K}}}_{S_{\si{K}}}(K_{\si{\infty}})$ 
est donné par:\smallskip

\centerline{$\rg_\Zl{}^\Gamma \C^{T_{\si{K}}}_{S_{\si{K}}}(K_{\si{\infty}}) = \deg_\ell \check T{\si{K}} + |\bar R_{\si{K^+}}^-|$.}
\end{Cor}

\Preuve Écrivant la pseudo-décomposition $\C^S_T(K_\%) \sim \Lambda^{\rho^{\si{S}}_{\si{T}}}\oplus \,\T^S_T(K_\%) $, nous avons immédiatement:\smallskip

\centerline{$^\Gamma \C^{T_{\si{K}}}_{S_{\si{K}}}(K_{\si{\infty}})
\sim \Zl^{\rho_{\si{S}}^{\si{T}}}\oplus \,^\Gamma \T^{T_{\si{K}}}_{S_{\si{K}}}(K_{\si{\infty}})
\sim \Zl^{\rho_{\si{S}}^{\si{T}}}\oplus \, \T^{T_{\si{K}}}_{S_{\si{K}}}(K_{\si{\infty}})^\Gamma
= \Zl^{\rho_{\si{S}}^{\si{T}}}\oplus \, \C^{T_{\si{K}}}_{S_{\si{K}}}(K_{\si{\infty}})^\Gamma$,}\smallskip

\noindent avec $\rho^{\si{S}}_{\si{T}}=\deg_\ell \check T_{\si{K}}$ et $\rg_\Zl\,\C_{\hat S_{\si{K}}}^{\check T_{\si{K}}}(K_\%)^\Gamma  = |\bar R_{\si{K^+}}^-|$.

\newpage
%%%%%%%%%%%%%%%%%%%%%%%%%%%%%%%%%%%%%%%%%%%%%%%%%%%%%%%%%%%%%%%%%%%%%%%
\section{Équivalence des conjectures faible et forte}
%%%%%%%%%%%%%%%%%%%%%%%%%%%%%%%%%%%%%%%%%%%%%%%%%%%%%%%%%%%%%%%%%%%%%%%

Supposons encore $K$ extension quadratique totalement imaginaire d'un sous-corps $K^+$ totalement réel; mais partons cette fois d'une partition arbitraire $S_{\si{K}} \sqcup T_{\si{K}}$ de $Pl_{\si{K}}^\ell$. Notons toujours $\tau$ la conjugaison complexe et posons: $\hat S_{\si{K}}=S_{\si{K}} \cup S_{\si{K}}^\tau$ et $\check S_{\si{K}}=S_{\si{K}} \cap S_{\si{K}}^\tau$ , comme dans le Théorème \ref{TP2}; et écrivons de même: $\hat T_{\si{K}}=T_{\si{K}} \cup T_{\si{K}}^\tau$ et $\check T_{\si{K}}=T_{\si{K}} \cap T_{\si{K}}^\tau$ en omettant l'indice $K$ dans ce qui suit.\par
Observons que $\hat S_{\si{K}} \sqcup \check T_{\si{K}}$ et  $\check S_{\si{K}} \sqcup \hat T_{\si{K}}$ forment des partitions de $P^\ell_{\si{K}}$ stables par la conjugaison $\tau$.\smallskip

Considérons alors la surjection canonique $f^{\si{S}}_{\si{T}}:\;\C^{\check S}_{\hat T}(K_\%)\to\,\C^{S}_{T}(K_\%)$. Son noyau $\Ker f^{\si{S}}_{\si{T}}$ est engendré conjointement par les sous-groupes de décomposition $\D_\l$ des places $\l$ de $S\setminus\check S$ et par les sous-groupes d'inertie $\I_\l$ des places $\l$ de $\hat T\setminus T$ dans $\,\C^{\check S}_{\hat T}(K_\%)$.
Regardons l'image par $f^{\si{S}}_{\si{T}}$ du sous-module réel $\,\C^{\check S}_{\hat T}(K_\%)^+=\,\C^{\check S}_{\hat T}(K_\%)^{1+\tau}\sim \,\C^{\check S}_{\hat T}(K_\%^+)$. Nous avons:\smallskip

\centerline{$\C^{S}_{T}(K_\%)/f^{\si{S}}_{\si{T}}(\C^{\check S}_{\hat T}(K_\%)^+)
\simeq\,\C^{\check S}_{\hat T}(K_\%)/\big(\Ker f^{\si{S}}_{\si{T}}\;\,\C^{\check S}_{\hat T}(K_\%)^+ \big)$.}\smallskip

Or, par construction, le quotient à droite est annulé par $1+\tau$. Ainsi, puisque les $\D_\l$ (pour $\l\in S\setminus\check S$) et les $\I_\l$ (pour $\l\in \hat T\setminus\ T$) y ont une image triviale, il en est de même de leurs conjugués respectifs $\D^\tau_\l=\D^\ph_{\l^\tau}$ et $\I^\tau_\l=\I^\ph_{\l^\tau}$. Et il vient donc:\smallskip

\centerline{$\C^{\check S}_{\hat T}(K_\%)/\big(\Ker f^{\si{S}}_{T}\,.\,\C^{\check S}_{\hat T}(K_\%)^+ \big)
\simeq \,\C^{\hat S}_{\check T}(K_\%)^-$,}\smallskip

\noindent où $\,\C^{\hat S}_{\check T}(K_\%)^-$ désigne le plus grand quotient de $\,\C^{\hat S}_{\check T}(K_\%)$ annulé par $1+\tau$.\smallskip

En résumé, nous avons la suite exacte courte:\smallskip

\centerline{$1 \longrightarrow \,\F^{\check S}_{\hat T}(K_\%)^+  \longrightarrow \,\C^{S}_{T}(K_\%) \longrightarrow  \,\C^{\hat S}_{\check T}(K_\%)^- \longrightarrow 1$,}\smallskip

\noindent avec $ \,\F^{\check S}_{\hat T}(K_\%)^+=  \,\C^{\check S}_{\hat T}(K_\%)^+/\big( \,\C^{\check S}_{\hat T}(K_\%)^+\cap \Ker f^{S}_{T}\big)$; puis, en prenant les points fixes par $\Gamma$:\smallskip

\centerline{$1 \longrightarrow \,\F^{\check S}_{\hat T}(K_\%)^{\Gamma\,+}  \longrightarrow \,\C^{S}_{T}(K_\%)^\Gamma \longrightarrow  \,\C^{\hat S}_{\check T}(K_\%)^{\Gamma\,-}$.}\smallskip

Supposons maintenant que $K$ satisfasse la conjecture cyclotomique {\em faible}, c'est à dire, comme établi dans \cite{J55}, à la fois la conjecture de Leopoldt et celle de Gross-Kuz'min pour $\ell$.
\begin{itemize}
\item À gauche, la conjecture de Leopoldt assure la finitude de $\,\C^{\check S}_{\hat T}(K_\%)^+$, donc de $\,\F^{\check S}_{\hat T}(K_\%)^{+}$; et finalement celle du sous-module des points fixes  $\,\F^{\check S}_{\hat T}(K_\%)^{\Gamma\,+}$.
\item À droite, la conjecture cyclotomique {\em  faible} (en fait la conjecture de Gross-Kuz'min) assure celle de $\,\C^{\hat S}_{\check T}(K_\%)^{\Gamma\,-}$.
\end{itemize}
\noindent Il suit de là que le groupe médian $\,\C^{S}_{T}(K_\%)^\Gamma$ est lui-même fini; autrement dit que $K$ vérifie la conjecture cyclotomique {\em forte} pour $\ell$. Ainsi:

\begin{Th}[\bf Équivalence des conjectures]\label{TP3}
Pour tout corps de nombres $K$ à conjugaison complexe et tout nombre premier $\ell$ fixé, les trois assertions suivantes sont équivalentes:
\begin{itemize}
\item $K$ vérifie les conjectures de Leopoldt et de Gross-Kuz'min.
\item $K$ vérifie la conjecture cyclotomique faible.
\item $K$ vérifie la conjecture cyclotomique forte.
\end{itemize}
\end{Th}

Comme vu dans les sections précédentes, la conjecture faible entraîne la validité de la formule des rangs pour les partitions $R \sqcup S \sqcup T$ de $Pl^\ell$ stables  par conjugaison; et la conjecture forte l'entraîne indépendamment de cette restriction. Inversement la formule écrite pour les seules partitions vérifiant $R=\emptyset$ exprime précisément ces mêmes conjectures.\smallskip

Avec les notations de Théorèmes \ref{TP1} et \ref{TP2}, il vient ainsi:

\begin{Sco}
Sont encore équivalentes aux précédentes chacune des deux assertions suivantes:
\begin{itemize}
\item Pour toute partition stable par conjugaison complexe $P_{\si{K}}^\ell=R_{\si{K}} \sqcup S_{\si{K}} \sqcup T_{\si{K}}$ de l'ensemble $Pl^\ell_{\si{K}}$, on a les identités de rang:\smallskip

\centerline{$\rg_{\Zl}\,\C^{T_{\si{K}}}_{S_{\si{K}}}(K^+_{\si{\infty}})^\Gamma =0$ \qquad \& \qquad $\rg_{\Zl}\,\C^{T_{\si{K}}}_{S_{\si{K}}}(K_{\si{\infty}})^\Gamma =|R_{\si{K^+}}^-|$.}\smallskip

\item Pour toute partition  $Pl_{\si{K}}^\ell=R_{\si{K}} \sqcup S_{\si{K}} \sqcup T_{\si{K}}$ de l'ensemble des places au-dessus de $\ell$, on a les identités de rang:

\centerline{$\rg_{\Zl}\,\C^{\check T_{\si{K}}}_{\hat S_{\si{K}}}(K^+_{\si{\infty}})^\Gamma =0$ \qquad \& \qquad $\rg_{\Zl}\,\C^{T_{\si{K}}}_{S_{\si{K}}}(K_{\si{\infty}})^\Gamma =|\bar R_{\si{K^+}}^-|$.}
\end{itemize}
\end{Sco}

\newpage
%%%%%%%%%%%%%%%%%%%%%%%%%%%%%%%%%%%%%%%%%%%%%%%%%%%%%%%%%%%%%%%%%%%%%%%
\section{Interprétation par le corps de classes $\ell$-adique}
%%%%%%%%%%%%%%%%%%%%%%%%%%%%%%%%%%%%%%%%%%%%%%%%%%%%%%%%%%%%%%%%%%%%%%%

Regardons maintenant le groupe $^\Gamma\C^{T_{\si{K}}}_{S_{\si{K}}}(K_{\si{\infty}})$ par la théorie $\ell$-adique du corps de classes, telle qu'exposée dans \cite{J31}: le nombre premier $\ell$ et le corps $K$ étant supposés fixés, pour chaque place non complexe $\p$ de $K$ notons 
$\,\U^{\phantom{*}}_\p$ et $\,\wU^{\phantom{*}}_\p$ les sous-groupes unités, respectivement au sens habituel et logarithmique, du $\ell$-adifié $\R_\p=\varprojlim \,K_\p^{\si{\times}}/K_\p^{\si{\times\ell^n}}$ du groupe multiplicatif du complété $K_\p$. 
Écrivons de même $\J=\prod_\p^{\si{\rm res}}\R_\p$ le $\ell$-adifié du groupe des idèles;  $\,\U=\prod_\p\U_\p$ et $\,\wU=\prod_\p\wU_\p$ ses sous-groupes unités (au sens habituel et logarithmique); et $\R=\Zl\otimes_\ZZ K^{\si{\times}}$ le sous-groupe principal de $\J$.\smallskip

Pour tout $P\subset Pl_K$, notons  $\R_P=\prod_{\p\in P}\R_\p$; puis $\,\U_P=\prod_{\p\in P}\U_\p$ et $\,\U^P=\prod_{\p\notin P}\U_\p$. Écrivons de même $\,\wU_P=\prod_{\p\in P}\wU_\p$ et $\,\wU^P=\prod_{\p\notin P}\wU_\p$. Soient enfin $\D_P=\R_P/\U_P$ et $\wD_P=\R_P/\wU_P$ les groupes de diviseurs (respectivement aux sens habituel et logarithmique). Cela étant, il vient:

\begin{Th}\label{TP4}
Soient $K$ un corps de nombres et $R \sqcup S \sqcup T$ une partition de l'ensemble $L$ des places de $K$ au-dessus de $\ell$. Avec les notations précédentes le groupe de Galois $\G^T_S=\Gal(H^T_S(K_\%/K)/K)$ de la pro-$\ell$-extension $S$-décomposée $T$-ramifiée maximale de $K_\%$ qui est abélienne sur $K$ est donné à un fini près par le pseudo-isomorphisme:\smallskip

\centerline{$\G_S^T \sim \wD_R\wD_S\,\U_T/(\wt\nu_{\si{RS}}\times p_{\si{T}})(\E_S)$,}\smallskip

\noindent où $\,\E_S$ est le $\ell$-adifié du groupe des $S$-unités de $K$;
$\wt\nu_{\si{RS}}=(\wi\nu_\l)_{\l\in R\cup S}$ est la famille des valuations logarithmiques aux places de $R \sqcup S$; et $p_{\si{T}}$ est le morphisme de semi-localisation aux places de $T$.
\end{Th}

\Preuve
Rappelons que  $\,\wU_\p$ est le groupe de normes attaché à la $\Zl$-extension cyclotomique de $K_\p$ et qu'on a $\,\U_\p=\,\wU_\p$, pour $\p\nmid \ell$. Posons $\,\bar\U_\p=\,\U_\p \cap\,\wU_\p$. Comme $\R_S\,\U^S\R=\R_S\,\U_R\,\U_T\,\U^L\R$ est d'indice fini dans $\J$, puisque $\J/\R_S\,\U^S\R$ s'identifie au $\ell$-groupe des $S$-classes d'idéaux, il vient:\smallskip

\centerline{$\G_S^T = \J / \,\wU_S\,\bar\U^T\R \sim  \R_S\,\U_R\,\U_{\,T}\,\U^L\R  /\,\wU_S\,\bar\U_R\,\U^L\R 
\simeq \R_S\,\U_R\,\U_{\,T} /\,\wU_S\,\bar\U_R(\R_S\,\U_R\,\U_{\,T}\cap \,\U^L\R)$;}\smallskip

\noindent ce qui donne la formule annoncée, puisque les idèles principaux qui interviennent au dénominateur sont les $S$-unités, et que l'on a, par ailleurs: $\R_S/\,\wU_S\simeq\wD_S$ et $\,\U_R/\,\bar\U_R \sim \R_R/\,\wU_R \simeq \wD_R$.

\begin{Cor}\label{CF}
Lorsque $K$ est extension quadratique totalement imaginaire d'un sous-corps $K^+$ totalement réel, sous les conjectures de Leopoldt et de Gross-Kuz'min, le pseudo-isomorphisme\smallskip

\centerline{$^\Gamma\C^{T_{\si{K}}}_{S_{\si{K}}}(K_{\si{\infty}}) \sim \,\G^{\check T_{\si{K}}\,-}_{\hat S_{\si{K}}}
\sim \,\wD_{\bar R_{\si{K}}}^- \times \,\U_{\,\check T_{\si{K}}}^-$ \quad 
donne directement: \quad $\rg_\Zl{}^\Gamma \C^{T_{\si{K}}}_{S_{\si{K}}}(K_{\si{\infty}}) = \deg_\ell \check T_{\si{K^+}} + |\bar R_{\si{K^+}}^-|$,}\smallskip

\noindent avec $\bar R_{\si{K}}= L_{\si{K}}\setminus(\check T_{\si{K}} \sqcup \hat S_{\si{K}})$ pour toute partition $R_{\si{K}} \sqcup S_{\si{K}} \sqcup T_{\si{K}}$ de $L_{\si{K}}=Pl_{\si{K}}^\ell$.

\end{Cor}

\Preuve 
Reprenons le schéma de démonstration du Théorème \ref{TP3} en partant cette fois de la surjection canonique $g^{T_{\si{K}}}_{S_{\si{K}}}: \,\G^{\check T_{\si{K}}}_{\hat S_{\si{K}}} \to \,\G^{T_{\si{K}}}_{S_{\si{K}}}$. Nous obtenons alors la suite pseudo-exacte courte:\smallskip

\centerline{$1 \longrightarrow \big(\G_{\check S_{\si{K}}}^{\hat T_{\si{K}}}\big)^+/ \big(\G_{\check S_{\si{K}}}^{\hat T_{\si{K}}}\cap \Ker g_{S_{\si{K}}}^{T_{\si{K}}}\big)^+  \longrightarrow \,\G_{S_{\si{K}}}^{T_{\si{K}}} \longrightarrow  \big(\G_{\hat S_{\si{K}}}^{\check T_{\si{K}}}\big)^- \longrightarrow 1$,}\smallskip

\noindent où apparaissent les composantes réelles ou imaginaires des modules de droite et de gauche, ce qui permet de définir partie réelle et partie imaginaire du groupe $\,\G_{S_{\si{K}}}^{T_{\si{K}}}$ alors même qu'il n'est pas {\em a priori} stable par conjugaison.\smallskip
\begin{itemize}
\item Côté réel, la conjecture de Leopoldt nous donne immédiatement: $\big(\G_{\check S_{\si{K}}}^{\hat T_{\si{K}}}\big)^+\sim \Gamma \simeq \Zl$; i.e. $^\Gamma\C^{T_{\si{K}}}_{S_{\si{K}}}(K_{\si{\infty}})^+\sim 1$; et finalement: $^\Gamma\C^{T_{\si{K}}}_{S_{\si{K}}}(K_{\si{\infty}} ) \sim  \big(\G_{\hat S_{\si{K}}}^{\check T_{\si{K}}}\big)^-$.\smallskip

\item Côté imaginaire, il vient, d'après le Théorème \ref{TP4}: $\G^{\check T_{\si{K}}}_{\hat S_{\si{K}}}{}^- \!\sim \wD_{\bar R_{\si{K}}}^-\wD_{\hat S_{\si{K}}}^-\,\U_{\,\check T_{\si{K}}}^-/(\wt\nu_{\si{\bar R\hat S}}\times p_{\si{\check T}})(\E_{\hat S_{\si{K}}}^-)$. Or, sous la conjecture de Gross-Kuz'min, le $\Zl$-module des $\hat S_{\si{K}}$-unités imaginaires s'envoie pseudo-injectivement dans $\wD_{\hat S_{\si{K}}}^-$, puisque les unités logarithmiques sont réelles. Il suit donc: $\,\E_{\hat S_{\si{K}}}^-\sim\wD_{\hat S_{\si{K}}}^-$; puis  $\,\G^{\check T_{\si{K}}}_{\hat S_{\si{K}}}{}^- \!\sim \wD_{\bar R_{\si{K}}}^-\,\U_{\,\check T_{\si{K}}}^-$, comme annoncé.
\end{itemize}\smallskip 

On retrouve ainsi très simplement l'expression du $\Zl$-rang donnée par le Théorème \ref{TP2}, ce qui fournit une démonstration alternative du Théorème d'équivalence \ref{TP3}.

\newpage
%%%%%%%%%%%%%%%%%%%%%%%%%%%%%%%%%%%%%%%%%%%%%%%%%%%%%%%%%%%%%%%%%%%%%%%
%%%%%%%%%%%%%%%%%%%%%%%%%%%%%%%%%%%%%%%%%%%%%%%%%%%%%%%%%%%%%%%%%%%%%%%

\noindent {\em Commentaires bibliographiques}\smallskip

Il est formulé dans \cite{J10} une conjecture générale sur l'indépendance $\ell$-adique de nombres algébriques, vérifiée dans les corps abéliens, qui implique en particulier les conjectures de Leopoldt et de Gross-Kuz'min et en factorise les preuves transcendantes (partielles) classiques.

 Les $\ell$-groupes de $S$-classes $T$-infinitésimales sont étudiés dans \cite{J18} (Ch. II, \S 2). On y trouve notamment la suite exacte des classes ambiges évoquée dans la section 1.\par

Les $\ell$-groupes de classes logarithmiques ont été introduits dans \cite{J28}. Leur calcul est maintenant implanté dans {\sc pari} (cf. \cite{BJ}). L'interprétation logarithmique de la conjecture de Gross-Kuz'min est donnée dans \cite{J28} sous l'appellation initiale de conjecture de Gross généralisée.\par

Une formulation équivalente de la conjecture cyclotomique en termes de points fixes de $(S,T)$-modules d'Iwasawa a été avancée dans \cite{LS} par Lee et Seo. Lee et Yu en ont tiré dans \cite{LY} pour $\ell$ impair un calcul du rang analogue à celui donné ici. Leur démonstration indépendante est plus laborieuse en l'absence des simplifications apportées par l'introduction des unités logarithmiques.\par
 
Les principaux résultats de la Théorie $\ell$-adique du corps de classes introduite dans \cite{J18} sont présentés dans \cite{J31}. On peut aussi se reporter au livre de Gras \cite{Gra2}. La semi-simplicité du module d'Iwasawa standard dans une $\Zl$-extension cyclotomiquement ramifiée est discutée dans \cite{JSa}.\par

Enfin, le calcul des invariants d'Iwasawa attachés aux $\ell$-groupes de $S$-classes $T$-infinitésimales est développé dans \cite{J43,JMa,JMP} en liaison avec les identités de dualité de Gras. L'article \cite{JMa} contient une erreur, reproduite dans \cite{J43} mais corrigée dans \cite{JMP}, qui ne concerne heureusement que l'invariant $\lambda^S_T$. Elle est sans incidence sur les résultats présentés ici.

\def\refname{\normalsize{\sc  Références}}

\medskip\noindent
{\small
\begin{tabular}{l}
Institut de Mathématiques de Bordeaux \\
Université de {\sc Bordeaux} \& CNRS \\
351 cours de la libération\\
F-33405 {\sc Talence} Cedex\\
courriel : Jean-Francois.Jaulent@math.u-bordeaux.fr\\
\url{https://www.math.u-bordeaux.fr/~jjaulent/}
\end{tabular}
}

\end{document}